\documentclass[preprint,12pt]{elsarticle}

\usepackage{amssymb}
\usepackage{amsthm}

\newcommand{\sysn}{\left\{\begin{array}{rcl}}
\newcommand{\sysk}{\end{array}\right.}

\newcommand{\ra}{\rangle}
\newcommand{\la}{\langle}

\renewcommand{\le}{\leqslant}
\renewcommand{\ge}{\geqslant}

\newtheorem{theorem}{Theorem}[section]

\theoremstyle{example}
\newtheorem{example}[theorem]{Example}

\theoremstyle{definition}
\newtheorem{definition}[theorem]{Definition}
\newtheorem{remark}[theorem]{Remark}
\newtheorem{corollary}[theorem]{Corollary}

\begin{document}

\begin{frontmatter}



\title{The $C$-compact-open topology on function spaces}


\author{Alexander V. Osipov}

\ead{OAB@list.ru}

\address{Ural Federal
 University, Institute of Mathematics and Mechanics, Ural Branch of the Russian Academy of Sciences, 16,
 S.Kovalevskaja street, 620219, Ekaterinburg, Russia}

\begin{abstract}
This paper studies the $C$-compact-open topology on the set $C(X)$
of all real-valued continuous functions on a Tychonov space $X$
and compares this topology with several well-known and lesser
known topologies. We investigate the properties $C$-compact-open
topology on the set $C(X)$ such as submetrizable, metrizable,
separable and second countability.
\end{abstract}

\begin{keyword}
set-open topology \sep $C$-compact subset \sep compact-open
topology \sep topological group \sep submetrizable


\MSC 54C40 \sep 54C35 \sep 54D60 \sep 54H11 \sep 46E10

\end{keyword}

\end{frontmatter}



\section{Introduction}

 The set-open topology on a family $\lambda$ of nonempty subsets of the set
 $X$ (the $\lambda$-open topology) is a generalization of the compact-open topology
 and of the topology of pointwise convergence. This topology was first introduced
 by Arens and Dugundji~{\cite{ardug}}.

 All sets of the form $[F, U]=\{f\in C(X):\ f(F)\subseteq U\}$, where $F\in\lambda$ and $U$ is an open subset
 of real line $\mathbb R$, form a subbase of the $\lambda$-open topology.

 The topology of uniform convergence is given by a base at each point $f\in C(X)$.
 This base consists of all sets $\{ g\in C(X):\ \sup\limits_{x\in X}|g(x)-f(x)|<~\varepsilon\}$.
 The topology of uniform convergence on elements of a family $\lambda$
 (the $\lambda$-topology), where $\lambda$ is a fixed family of non-empty subsets of the
 set~$X$, is a natural generalization of this topology.
 All sets of the form $\{g\in C(X): \sup\limits_{x\in F}|g(x)-f(x)|<\varepsilon\}$,
  where $F\in\lambda$ and $\varepsilon >0$,
 form a base of the $\lambda$-topology at a point $f\in C(X)$.

Note that a $\lambda$-open topology coincides with a
$\lambda$-topology, when the family  $\lambda$ consists of all
finite (compact, countable compact, pseudocompact, sequentially
compact, $C$-compact) subsets of $X$. Therefore $C(X)$ with the
topology of pointwise convergence (compact-open, countable
compact-open, sequentially compact-open, pseudocompact-open,
$C$-compact-open topology) is a locally convex topological vector
space.

Moreover, if a $\lambda$-open topology coincides with a
$\lambda$-topology, then $\lambda$ consists of $C$-compact subsets
of space $X$ and the space $C_{\lambda}(X)$ is a topological
algebra under the usual operations of addition and multiplication
(and multiplication by scalars).

\section{Main definitions and notation}

 In this paper, we consider the space $C(X)$ of all real-valued continuous functions
 defined on a Tychonov space~$X$.
  We denote by $\lambda$ a family of non-empty subsets
 of the set~$X$. We use the following notation for various topological spaces
 with the underlying set $C(X)$:

 $C_{\lambda}(X)$ for the $\lambda$-open topology,

 $C_{\lambda, u}(X)$ for the $\lambda$-topology.

 The elements of the standard subbases of the $\lambda$-open topology and $\lambda$-topology
 will be denoted as follows:

 $[F,\,U]=\{f\in C(X):\ f(F)\subseteq U\}$,

 $\la f,\,F,\,\varepsilon\ra=\{g\in C(X):\ \sup\limits_{x\in F}|f(x)-g(x)|<\varepsilon\}$,
  where $F\in\lambda$, $U$ is an open subset  of $\mathbb R$ and $\varepsilon>0$.

 If $X$ and $Y$ are  any two topological spaces with the same underlying set, then we use the notation
 $X=Y$, $X\le Y$, and $X<Y$ to
 indicate, respectively, that $X$ and $Y$ have the same topology,
 that the topology on $Y$ is finer than or equal to the topology on
  $X$, and that the topology on $Y$ is strictly finer than the topology on $X$.

 The closure of a set $A$ will be denoted by $\overline{A}$; the symbol
 $\varnothing$ stands for the empty set. As usual, $f(A)$ and $f^{-1}(A)$ are the image and
 the complete preimage of the set $A$ under the mapping~$f$,
 respectively. The constant zero function defined on $X$ is
 denoted by $0$, more precisely by $0_{X}$. We call it the
 constant zero function in $C(X)$.

 We denote by $\mathbb{R}$ the real line
 with the natural topology.

  We recall that a subset of $X$ that is the
 complete preimage of zero for a certain function from~$C(X)$ is called a zero-set.
A subset $O$ of a space $X$ is called functionally open (or a
cozero-set) if $X\setminus O$ is a zero-set. A family $\lambda$ of
non-empty subsets of a topological space $(X,\tau)$ is called a
$\pi$-network for $X$ if for any nonempty open set $U\in\tau$
there exists $A\in \lambda$ such that $A\subset U$.

Let $\overline{\lambda}=\{\overline{A} : A\in \lambda \}$. Note
that the same set-open topology is obtained if $\lambda$ is
replaced by $\overline{\lambda}$. This is because for each $f\in
C(X)$ we have $f(\overline{A})\subseteq \overline{f(A)}$ and,
hence, $\overline{f(\overline{A})}=\overline{f(A)}$. Consequently,
$C_{\overline{\lambda}}(X)=C_{\lambda}(X)$. From now on, $\lambda$
denotes a family of non-empty closed subsets  of the set~$X$.

  Throughout this paper, a family $\lambda$ of nonempty subsets of the set $X$ is
  a $\pi$-network. This condition is equivalent to the space $C_{\lambda}(X)$ being a Hausdorff space.
 The set-open topology does not change when $\lambda$  is replaced with the finite
unions of its elements. Therefore we assume that $\lambda$ is
closed under finite unions of its elements.

Recall that a subset $A$ of a space $X$ is called
 {\it $C$-compact subset $X$} if, for any real-valued function~$f$ continuous on $X$,
 the set $f(A)$ is compact in ~${\mathbb{R}}$.

Note (see Theorem 3.9 in {\cite{os1}}) that the set $A$  is a
$C$-compact  subset of $X$ if and only if every countable
functionally open (in $X$) cover of $A$ has a finite subcover.

The remaining notation can be found in~{\cite{enge}}.

\section{Topological-algebraic properties of function spaces}

Interest in studying the $C$-compact topology generated by a
Theorem 3.3 in ~{\cite{os2}} which characterizes some
topological-algebraic properties of the set-open topology. It
turns out if $C_{\lambda}(X)$ is a paratopological group (TVS ,
locally convex TVS) then the family $\lambda$ consists of
$C$-compact subsets of $X$.

Given a family $\lambda$ of non-empty subsets of $X$, let
$\lambda(C)=\{A\in \lambda$ :  for every $C$-compact subset $B$ of
the space $X$ with $B\subset A$, the set $[B,U]$ is open in
$C_{\lambda}(X)$ for any open set $U$ of the space $\mathbb{R}
\}$.

 Let $\lambda_m$ denote the maximal with respect to inclusion
 family, provided that $C_{\lambda_m}(X)=C_{\lambda}(X)$. Note that a family
 $\lambda_m$ is unique for each family $\lambda$.

 A family $\lambda$ of $C$-compact subsets of $X$ is said to be
 hereditary with respect to $C$-compact subsets if it satisfies
 the following condition: whenever $A\in \lambda$ and $B$ is
 a $C$-compact (in $X$) subset of $A$, then $B\in \lambda$ also.

We look at the properties of the family $\lambda$ which imply that
the space $C_{\lambda}(X)$ with the set-open topology is a
topological algebra under the usual operations of addition and
multiplication (and multiplication by scalars).

The following theorem is a generalization of Theorem 3.3 in
~{\cite{os2}}.

\begin{theorem}\label{th1} For a space $X$, the following statements are
equivalent.

\begin{enumerate}

\item  $C_{\lambda}(X)=C_{\lambda, u}(X).$

\item  $C_{\lambda}(X)$ is a paratopological group.

\item  $C_{\lambda}(X)$ is a topological group.

\item  $C_{\lambda}(X)$ is a topological vector space.

\item $C_{\lambda}(X)$ is a locally convex topological vector
space.

\item $C_{\lambda}(X)$ is a topological ring.

\item $C_{\lambda}(X)$ is a topological algebra.

\item  $\lambda$ is a family of\, $C$-compact sets and
$\lambda=\lambda(C)$.

\item $\lambda_m$ is a family of\, $C$-compact sets and it is
hereditary with respect to $C$-compact subsets.

\end{enumerate}

\end{theorem}

\begin{proof}
Equivalence of the statements (1), (3), (4), (5) and (8) proved in
~{\cite[Theorem 3.3]{os2}}.

Note that in the proof of Lemmas 3.1 and 3.2 in ~{\cite{os2}} used
only the condition that the space $X$ is a paratopological space.
Thus (2)$\Rightarrow$(8).

(8)$\Rightarrow$(7). As (8)$\Leftrightarrow$(4), we only need to
show that continuous operation of multiplication. Really, let
$\beta$ be a neighborhood filter of zero function in $C(X)$. Let
$W=[A, V]\in \beta$, where $A\in \lambda$ and $V$ is an open set
of the space $\mathbb{R}$. Then there is an open set $V^1$ such
that $V^1*V^1\subset V$. Show that $W^1=[A, V^1]$ such that
$W^1*W^1\subset W$. Indeed $W^1*W^1=\{f*g : f\in W^1, g\in W^1
\}=\{ f*g : f(A)\subset V^1$ and $g(A)\subset V^1\}$. Clearly that
$f(x)*g(x)\in V^1*V^1$ for each $x\in A$. Therefore
$(f*g)(A)\subset V$ and $W^1*W^1\subset W$.

It remains to prove that if $W=[A, V]\in \beta$ and $f\in C(X)$
then there is an open set $V^1\ni 0$ such that $f(A)*V^1\subset V$
and $V^1*f(A)\subset V$. Indeed let $g=f*h$ and $g^1=h^1*f$ where
$h, h^1\in W^1$. Then $g(x)=f(x)*h(x)\in f(A)*V^1$ and
$g^1(x)=h^1(x)*f(x)\in V^1*f(A)$ for each $x\in A$. Note that
$g(A)\subset V$ and $g^1(A)\subset V$.

(8)$\Rightarrow$(9). Since $C_{\lambda_m}(X)=C_{\lambda}(X)$ then
$C_{\lambda_m}(X)$ is a topological group and $\lambda_m$ is a
family of\, $C$-compact sets and consequently,
$\lambda_{m}=\lambda_{m}(C)$. But if the set $[B,U]$ is open in
$C_{\lambda_m}(X)$ for any open set $U$ of the space $\mathbb{R}$
then $B\in \lambda_m$.

Remaining implications is obviously and follow from Theorem 3.3
in~{\cite{os2}} and the definitions.

\end{proof}

\section{Comparison of topologies}

 In this section, we compare the $C$-compact-open topology with
 several well-known and lesser known topologies.

 We use the following notations to denote the particular families
 of $C$-compact subsets of $X$.

  $F(X)$ --- the collection of all finite subsets of $X$.

$MC(X)$ --- the collection of all metrizable compact subsets of
 $X$.

  $K(X)$ --- the collection of all compact subsets of $X$.

 $SC(X)$ --- the collection of all sequential compact subsets of $X$.

 $CC(X)$ --- the collection of all countable compact subsets of $X$.

 $PS(X)$ --- the collection of all pseudocompact subsets of $X$.

 $RC(X)$ --- the collection of all $C$-compact subsets of $X$.

Note that $F(X)\subseteq MC(X)\subseteq K(X)\subseteq
CC(X)\subseteq PS(X)\subseteq RC(X)$ and $MC(X)\subseteq
SC(X)\subseteq CC(X)$.
  When $\lambda=F(X)$, $MC(X)$, $K(X)$, $SC(X)$, $CC(X)$, $PS(X)$ or
 $RC(X)$, we call the corresponding $\lambda$-open topologies on
 $C(X)$ point-open, metrizable compact-open, compact-open,
 sequential compact-open, countable compact-open, pseudocompact-open and
 $C$-compact-open respectively. The corresponding spaces are
 denoted by $C_{p}(X)$, $C_{k}(X)$, $C_{c}(X)$,
 $C_{sc}(X)$, $C_{cc}(X)$, $C_{ps}(X)$ and $C_{rc}(X)$ respectively.

For the $C$-compact-open topology on $C(X)$, we take as subbase,
the family $\{[A,\,V] : A\in RC(X), V$ is open in $\mathbb R \}$;
and we denote the corresponding space by $C_{rc}(X)$.

We obtain from Theorem \ref{th1}  the following result.

\begin{theorem} For any space $X$ and $\lambda\in \{F(X)$, $MC(X)$, $K(X)$, $SC(X)$, $CC(X)$, $PS(X)$, $RC(X) \}$, the $\lambda$-open topology
on $C(X)$ is same as the topology of uniform convergence on
elements of a family $\lambda$, that is,
$C_{\lambda}(X)=C_{\lambda, u}(X)$. Moreover, $C_{\lambda}(X)$ is
a Hausdorff locally convex topological vector space (TVS).
\end{theorem}

When $C(X)$ is equipped with the topology of uniform convergence
on $X$, we denote the corresponding space by $C_{u}(X)$.

\begin{theorem} For any space $X$,

$C_{p}(X)\leq C_{k}(X)\leq C_{c}(X)\leq C_{cc}(X)\leq
C_{ps}(X)\leq C_{rc}(X)\leq C_{u}(X)$

and

$C_{k}(X)(X)\leq C_{sc}(X)\leq C_{cc}(X)$.

\end{theorem}

Now we determine when these inequalities are equalities and give
examples to illustrate the differences.

\begin{example} \label{ex1}

Let $X$ be the set of all countable ordinals  $\{\alpha : \alpha <
\omega_1\}$  equipped with the order topology. The space $X$ is
sequential compact  and collectionwise normal, but not compact.
For this space $X$, we have  $C_{sc}(X)>C_{c}(X)$.

 Really, let $f=f_{0}$ and $U=(-1,1)$. Consider the
neighborhood $[X,U]$ of  $f$. Assume that there are a family of
neighborhoods $\{[A_i,U_i]\}_{i=1}^{n}$, where $A_i$
--- compact, and $f\in\bigcap\limits_{i=1}^n [A_i,U_i] \subset
[X,U]$. Then $\exists\alpha < \omega_1$ such that $\forall
\beta\in A_i\ \beta<\alpha$. Define function $g$: $g(\beta)=0$ for
$\beta\le\alpha$ and $g(\beta)=1$ for $\beta>\alpha$. Then $g\in
\bigcap\limits_{i=1}^n [A_i,U_i]$, but $g\notin [X,U]$, a
contradiction.

Note that for this space $X$, we have:

 $C_{p}(X)<C_{k}(X)=C_{c}(X)<C_{sc}(X)=C_{cc}(X)=C_{ps}(X)=C_{rc}(X)=C_u(X)$.

\end{example}

\begin{example} \label{ex2}

 Let $Y=\beta \mathbb N$ be Stone-C$\check{e}$ch compactification
 of natural numbers $\mathbb N$. Note that every sequential compact
 subset of $\beta \mathbb N$ is finite. For this space $Y$, we
 have:

$C_{p}(Y)=C_{k}(Y)=C_{sc}(Y)<C_{c}(Y)=C_{cc}(Y)=C_{ps}(Y)=C_{rc}(Y)=C_{u}(Y)$.

\end{example}

\begin{example} \label{ex3}
 Let  $Z=X\oplus Y$ where $X$ is the space of Example~\ref{ex1}
 and $Y$ is the space of Example~\ref{ex2}. Then the sequential compact-open topology is incomparable with
 the compact-open topology on the space $C(Z)$.
\end{example}

\begin{example} \label{ex4}
 Let  $X=I^{\mathfrak c}$  be the Tychonoff cube of weight $\mathfrak
c$. A space $X$ is compact and contains a dense sequential compact
subset. Thus, we have:

 $C_{p}(X)<C_{k}(X)<C_{c}(X)=C_{sc}(X)=C_{cc}(X)=C_{ps}(X)=C_{rc}(X)=C_{u}(X)$.

\end{example}

\begin{example} \label{ex5}
Let $X=\omega_1+1$ be the set of all ordinals $\leqslant \omega_1$
equipped with the order topology. The space $X$ is compact and
sequentially compact but not metrizable. Then, for space $X$ we
have:

$C_p(X)<C_k(X)<C_{sc}(X)=C_c(X)=C_{cc}(X)=C_{ps}(X)=C_{rc}(X)=C_u(X)$.

\end{example}

The following example is an example of the space in which every
sequentially compact and compact subset is finite.

\begin{example} \label{ex6}
 Let $K_0 = \mathbb N$. By using transfinite induction,
 we construct subspace of  $\beta \mathbb N$. Suppose that $K_\beta \subset
\beta \mathbb N$ is defined for each  $\beta < \alpha$ and
$|K_\beta|\leqslant \mathfrak{c}$. Then for each $A\in
[\bigcup_{\beta<\alpha} K_\beta]^\omega$ choose $x_A$ such that
$x_A$ is an accumulation point of the set $A$ in the space $\beta
\mathbb N$. Let $K_\alpha = \bigcup_{\beta<\alpha} K_\beta \cup
\{x_A\colon A\in [\bigcup_{\beta<\alpha} K_\beta]^\omega\}$. The
space $M=\bigcup_{\alpha<\omega_1} K_\alpha$ is countable compact
space in which every sequentially compact and compact subset is a
finite. Thus, we have:

$C_{p}(M)=C_{k}(M)=C_{sc}(M)=C_{c}(M)<C_{cc}(M)=C_{ps}(M)=C_{rc}(M)=C_{u}(M)$.
\end{example}

\begin{example} \label{ex7}
Let $\mathcal M$ be a maximal infinite family of infinite subsets
of $\mathbb N$ such that the intersection of any two members of
$\mathcal M$ is finite, and let $\Psi=\mathbb N\bigcup \mathcal
M$, where a subset $U$ of $\Psi$ is defined to be open provided
that for any set $M\in \mathcal M$, if $M\in U$ then there is a
finite subset $F$ of $M$ such that $\{M\}\bigcup M\setminus
F\subset U$. The space $\Psi$ is then a first-countable
pseudocompact Tychonov space that is not countably compact. The
space $\Psi$  is due independently to J. Isbell and S.
Mr$\acute{o}$wka.

Every compact, sequentially compact, countable compact subsets of
$\Psi$ has the form  $\bigcup_{i=1}^n (\{x_i\} \cup (x_i\setminus
S_i)) \cup S$, where $x_i\in E$, $|S_i| < \omega$, $|S| < \omega$.
Thus obtain the following relations:

$C_{p}(\Psi)<C_{k}(\Psi)=C_{sc}(\Psi)=C_{c}(\Psi)=C_{cc}(\Psi)<C_{ps}(\Psi)=C_{rc}(\Psi)=C_{u}(\Psi)$.
\end{example}

\begin{example} \label{ex8}

 Let  $X = \beta \mathbb N \oplus (\omega_1+1) \oplus M \oplus
\Psi$, where $M$ is the space of Example~\ref{ex6}
 and $\Psi$ is the space of Example~\ref{ex7}.
We have the following relations:

 $C_{p}(X)<C_{k}(X)<C_{sc}(X)<C_{c}(X)<C_{cc}(X)<C_{ps}(X)=C_{rc}(X)=C_{u}(X)$.

 \end{example}

\begin{example} \label{ex9}

 Let $G=(\omega_1) \oplus M \oplus \Psi \oplus I^\mathfrak{c}$,
where $M$ is the space of Example~\ref{ex6}
 and $\Psi$ is the space of Example~\ref{ex7}.
We have the following relations:

$C_{p}(G)<C_{k}(G)<C_{c}(G)<C_{sc}(G)<C_{cc}(G)<C_{ps}(G)=C_{rc}(G)=C_{u}(G)$.

\end{example}

\begin{example} \label{ex10}  Let $Y=[0,\omega_2]\times[0,\omega_1]\setminus\{(\omega_2,\omega_1)\}$,
with the topology $\tau$ generated by declaring open each point of
$[0,\omega_2)\times[0,\omega_1)$, together with the sets
$U_{P}(\beta)=\{(\beta,\gamma):
   \gamma\in ([0,\omega_1]\setminus P)$, where $P$ is finite and $(\beta,\omega_1)\notin P \}$
    and $V_{\alpha}(\beta)=\{(\gamma,\beta): \alpha<\gamma\leqslant\omega_2\}$.

Let $A=\{(\omega_2,\gamma): 0\leqslant \gamma<\omega_1\}$ and
$f\in C(Y)$. Suppose that $f(A)$ is not a closed set, then there
are $c\in \overline{f(A)}\setminus f(A)$ and sequence
$\{a_{n}\}\subset A$ such that $\{f(a_{n})\}\rightarrow c$. Since
$a_{n}=(\omega_2,\gamma_n)$, there is $\alpha_n$ such that
$f(\alpha,\gamma_n)=f(a_{n})$ for each $\alpha>\alpha_n$.
Moreover, there is $\beta$, such that
$f(\alpha,\gamma)=f(\omega_2,\gamma)$ for each $\gamma\in
[0,\omega_1]$ and $\alpha\geq\beta$. Clearly that
$f(\beta,\omega_1)=c$. Then there is $\delta$, such that
$f(\beta,\gamma)=c$  for each $\gamma\geq\delta$. It follows that
$f(\omega_2,\gamma)=c$, but $(\omega_2,\gamma)\in A$ and $c\notin
f(A)$, a contradiction. Thus, set $A$ is a $C$-compact subset of
the space $Y$.

Let $B$ be a nonempty pseudocompact subset of $Y$.
  Since $\{\alpha\}\times [0,\omega_1]$ is a clopen set (functionally
  open) for each $\alpha<\omega_2$, $([0,\omega_2]\times \{\beta\})\bigcap
  B$ has at most a finite number of points for each $\beta\leqslant\omega_1$.
It follows that $B$ is a compact subset of $Y$.

 As $A$ is the infinite set and closed and the
 pseudocompact subsets of $Y$ are compact and have at most a
 finite intersection with $A$, $A$ provides an example of a
 $C$-compact subset which is not contained in any closed
 pseudocompact subset of $Y$. Since $Y$ has infinite compact
 subsets, for this space we have

 $C_{c}(Y)=C_{ps}(Y)<C_{rc}(Y)$.

\end{example}

\begin{example} \label{ex11} Let $Z=Y\bigoplus G$,
where $Y$ is the space of Example~\ref{ex10}
 and $G$ is the space of Example~\ref{ex9}.
We have the following relations:

$C_{p}(Z)<C_{k}(Z)<C_{c}(Z)<C_{sc}(Z)<C_{cc}(Z)<C_{ps}(Z)<C_{rc}(Z)$.

\end{example}

Recall that a space $X$ is called submetrizable if $X$ admits a
weaker metrizable topology.

Note that for a subset $A$ in a submetrizable space $X$, the
following are equivalent:

(1)   $A$ is metrizable compact,

(2)   $A$ is compact,

(3)   $A$ is sequential compact,

(4)   $A$ is countable compact,

(3)   $A$ is pseudocompact,

(4)   $A$ is $C$-compact subset of $X$.

\begin{theorem}\label{th10}
 Let $X$ be a submetrizable space, then

$C_{k}(X)=C_{c}(X)=C_{sc}(X)=C_{cc}(X)=C_{ps}(X)=C_{rc}(X)$.

\end{theorem}

Similarly to Corollary 3.7 in {\cite{kund}} on the bounded-open
topology we have

\begin{theorem} For every space $X$,

\begin{enumerate}

\item $C_{c}(X)=C_{rc}(X)$ iff every closed $C$-compact subset of
$X$ is compact.

\item $C_{rc}(X)=C_{u}(X)$ iff $X$ is pseudocompact.

\end{enumerate}

\end{theorem}

\begin{proof}
(1) Note that for a subset $A$ of $X$, $\la f, \overline{A},
\epsilon \ra \subseteq \la f, A, \varepsilon \ra$. So if every
closed $C$-compact subset of $X$ is compact, then $C_{rc}(X)\leq
C_{c}(X)$. Consequently, in this case, $C_{rc}(X)=C_{c}(X)$.

Conversely, suppose that $C_{c}(X)=C_{rc}(X)$ and let $A$ be any
closed $C$-compact subset of $X$. So $\la 0, A, 1 \ra$ is open in
$C_{c}(X)$ and consequently, there exist a compact subset $K$ of
$X$ and $\varepsilon>0$ such that $\la 0, K, \varepsilon \ra
\subseteq
 \la 0, A, 1\ra$. If possible, let $x\in A\setminus K$. Then there exists a continuous
function $g : X \mapsto [0, 1]$ such that $g(x) = 1$ and $g(y) =
0$ $\forall y\in K$. Note that $g\in \la 0, K, \varepsilon \ra
\setminus \la 0, A, 1 \ra$ and we arrive at a contradiction.
Hence, $A\subseteq K$ and consequently, $A$ is compact.

(2) First, suppose that $X$ is pseudocompact. So for each $f\in
C(X)$ and each $\varepsilon > 0$, $\la f, X, \varepsilon \ra$ is a
basic open set in $C_{rc}(X)$ and consequently,
$C_{u}(X)=C_{rc}(X)$.

Now let $C_{rc}(X)=C_{u}(X)$. Since $\la 0, X, 1 \ra$ is a basic
neighborhood of the constant zero function $0$ in $C_{u}(X)$,
there exist a $C$-compact subset $A$ of $X$ and $\varepsilon> 0$
such that $\la 0, A, \varepsilon \ra \subseteq \la 0, X, 1 \ra$.
As before, by using the complete regularity of $X$, it can be
shown that we must have $X =\overline{A}$. But the closure of a
$C$-compact set is also $C$-compact set. Hence, $X$ is
pseudocompact.
\end{proof}

Note that for a closed subset $A$ in a normal Hausdorff space $X$,
the following are equivalent:

(1)   $A$ is countable compact,

(2)   $A$ is pseudocompact,

(3)   $A$ is $C$-compact subset of $X$.

\begin{corollary}
For any normal Hausdorff space $X$, $C_{c}(X)=C_{rc}(X)$ iff every
closed countable compact subset of $X$ is compact.
\end{corollary}

\section{Submetrizable and metrizable}

One of the most useful tools in function spaces is the following
concept of induced map. If $f : X \mapsto Y$ is a continuous map,
then the induced map of $f$, denoted by $f^* : C(Y)\mapsto C(X)$
is defined by $f^*(g)=g\circ f$ for all $g\in C(Y)$.

Recall that a map $f: X\mapsto Y$ , where $X$ is any nonempty set
and $Y$ is a topological space, is called almost onto if $f(X)$ is
dense in $Y$.

\begin{theorem}\label{th5} Let $f: X \mapsto Y$ be a continuous map between two
spaces $X$ and $Y$. Then

(1) $f^* : C_{rc}(Y) \mapsto C_{rc}(X)$ is continuous;

(2) $f^* : C(Y)\mapsto C(X)$ is one-to-one if and only if $f$ is
almost onto;

(3) if $f^* : C(Y)\mapsto C_{rc}(X)$ is almost onto, then $f$ is
one-to-one.
\end{theorem}

\begin{proof} : (1) Suppose $g\in C_{rc}(Y)$. Let $\la f^*(g), A,
\varepsilon \ra$ be a basic neighborhood of $f^*(g)$ in
$C_{rc}(X)$. Then $f^*(\la g, f(A), \varepsilon \ra)\subseteq \la
f^*(g), A, \epsilon \ra$ and consequently, $f^*$ is continuous.
(2) and (3) See Theorem 2.2.6 in ~{\cite{mcnt}}.

\end{proof}

\begin{remark}
(1) If a space $X$ has a $G_{\delta}$-diagonal, that is, if the
set $\{(x, x): x\in X\}$ is a $G_{\delta}$-set in the product
space $X\times X$, then every point in $X$ is a $G_{\delta}$-set.
Note that every metrizable space has a zero-set diagonal.
Consequently, every submetrizable space has also a
zero-set-diagonal.

(2) Every compact set in a submetrizable space is a
$G_{\delta}$-set. A space $X$ is called an $E_0$-space if every
point in the space is a $G_{\delta}$-set. So the submetrizable
spaces are $E_0$-spaces.
\end{remark}

For our next result, we need the following definitions.

\begin{definition} A completely regular Hausdorff space $X$ is called
$\sigma$-$C$-compact if there exists a sequence $\{A_{n}\}$ of
$C$-compact sets in $X$ such that $X=\bigcup_{n=1}^{\infty}
A_{n}$. A space $X$ is said to be almost $\sigma$-$C$-compact if
it has a dense $\sigma$-$C$-compact subset.

\end{definition}

\begin{theorem} For any space $X$, the following are equivalent.

\begin{enumerate}

\item  $C_{rc}(X)$ is submetrizable.

\item Every $C$-compact subset of $C_{rc}(X)$ is a
$G_{\delta}$-set in $C_{rc}(X)$.

\item Every compact subset of $C_{rc}(X)$ is a $G_{\delta}$-set in
$C_{rc}(X)$.

\item $C_{rc}(X)$ is an $E_0$-space.

\item $X$ is almost $\sigma$-$C$-compact.

\item $C_{rc}(X)$ has a zero-set-diagonal.

\item $C_{rc}(X)$ has a $G_{\delta}$-diagonal.

\end{enumerate}

\end{theorem}

\begin{proof} $(1)\Rightarrow(2)\Rightarrow(3)\Rightarrow(4)$ are all immediate.

$(4)\Rightarrow(5)$. If $C_{rc}(X)$ is an $E_0$-space, then the
constant zero function $0$ defined on $X$ is a $G_{\delta}$-set.
Let $\{0\}=\bigcap_{n=1}^{\infty} \la 0, A_{n}, \varepsilon \ra$
where each $A_{n}$ is $C$-compact subset in $X$ and
$\varepsilon>0$. We claim that $X=\overline{\bigcup_{n=1}^{\infty}
A_{n}}$.

Suppose that $x_{0}\in X\setminus \overline{\bigcup_{n=1}^{\infty}
A_{n}}$. So there exists a continuous function $f : X \mapsto [0,
1]$ such that $f(x)=0$ for all $x\in
\overline{\bigcup_{n=1}^{\infty} A_{n}}$ and $f(x_{0})=1$. Since
$f(x)=0$ for all $x\in A_n$, $f\in \la 0, A_n, \varepsilon \ra$
for all $n$ and hence, $f\in \bigcap_{n=1}^{\infty} \la 0, A_{n},
\varepsilon \ra=\{0\}$. This means $f(x) = 0$ for all $x\in X$.
But $f(x_0) = 1$. Because of this contradiction, we conclude that
X is almost $\sigma$-$C$-compact.

$(5)\Rightarrow(1)$. By Theorem 4.10 in~{\cite{kuga}} and Theorem
~\ref{th5}.

By Remark $(1)\Rightarrow(6)\Rightarrow(7)\Rightarrow(4)$.

\end{proof}

\begin{corollary}

Suppose that $X$ is almost $\sigma$-$C$-compact. If $K$ is a
subset of $C_{rc}(X)$, then the following are equivalent.

\begin{enumerate}

\item  $K$ is metrizable compact.

\item  $K$ is compact.

\item  $K$ is sequentially compact.

\item  $K$ is countable compact.

\item  $K$ is pseudocompact.

\item  $K$ is $C$-compact subset of $C_{rc}(X)$.

\end{enumerate}

\end{corollary}

 A space $X$ is said to be of (pointwise) countable type if each
(point) compact set is contained in a compact set having countable
character.

A space $X$ is a $q$-space if for each point $x\in X$, there
exists a sequence $\{U_{n} : n\in \mathbb N \}$ of neighborhoods
of $x$ such that if $x_{n}\in U_{n}$ for each $n$, then $\{x_{n} :
n\in \mathbb N\}$ has a cluster point. Another property stronger
than being a $q$-space is that of being an $M$-space, which can be
characterized as a space that can be mapped onto a metric space by
a quasi-perfect map (a continuous closed map in which inverse
images of points are countably compact). Both a space of pointwise
countable type and an $M$-space are $q$-spaces.

\begin{theorem} For any space $X$, the following are equivalent.

\begin{enumerate}

\item  $C_{rc}(X)$ is metrizable.

\item $C_{rc}(X)$ is of first countable.

\item $C_{rc}(X)$ is of countable type.

\item $C_{rc}(X)$ is of pointwise countable type.

\item $C_{rc}(X)$ has a dense subspace of pointwise countable
type.

\item $C_{rc}(X)$ is an $M$-space.

\item $C_{rc}(X)$ is a $q$-space.

\item $X$ is hemi-$C$-compact; that is, there exists a sequence of
$C$-compact sets $\{A_{n}\}$ in $X$ such that for any $C$-compact
subset $A$ of $X$, $A\subseteq A_{n}$ holds for some $n$.

\end{enumerate}

\end{theorem}

\begin{proof}
From the earlier discussions, we have
$(1)\Rightarrow(3)\Rightarrow(4)\Rightarrow(7)$,
$(1)\Rightarrow(6)\Rightarrow(7)$, and
$(1)\Rightarrow(2)\Rightarrow(7)$.

$(4)\Leftrightarrow(5)$. It can be easily verified that if $D$ is
a dense subset of a space $X$ and $A$ is a compact subset of $D$,
then $A$ has countable character in $D$ if and only if $A$ is of
countable character in $X$. Now since $C_{rc}(X)$ is a locally
convex space, it is homogeneous. If we combine this fact with the
previous observation, we have $(4)\Leftrightarrow(5)$.

$(7)\Rightarrow(8)$. Suppose that $C_{rc}(X)$ is a $q$-space.
Hence, there exists a sequence $\{U_{n} : n\in \mathbb N \}$ of
neighborhoods of the zero-function $0$ in $C_{rc}(X)$ such that if
$f_{n}\in U_{n}$ for each $n$, then $\{f_{n} : n\in \mathbb N \}$
has a cluster point in $C_{rc}(X)$. Now for each $n$, there exists
a closed $C$-compact subset $A_{n}$ of $X$ and $\epsilon_{n}>0$
such that $0\in \la 0, A_{n}, \epsilon_{n} \ra \subseteq U_{n}$.

Let $A$ be a $C$-compact subset of $X$. If possible, suppose that
$A$ is not a subset of $A_{n}$ for any $n\in \mathbb N$. Then for
each $n\in \mathbb N$, there exists $a_{n}\in A\setminus A_{n}$.
So for each $n\in \mathbb N$, there exists a continuous function
$f_{n} : X \mapsto [0, 1]$ such that $f_{n}(a_{n})=n$ and
$f_{n}(x)=0$ for all $x\in A_n$. It is clear that $f_n\in \la 0,
A_n, \epsilon_n \ra$. But the sequence $\{f_{n}\}_{n\in \mathbb
N}$ does not have a cluster point in $C_{rc}(X)$. If possible,
suppose that this sequence has a cluster point $f$ in $C_{rc}(X)$.
Then for each $k\in \mathbb N$, there exists a positive integer
$n_k>k$ such that $f_{n_{k}}\in \la f, A, 1 \ra$. So for all $k\in
N$, $f(a_{n_{k}})> f_{n_k} (a_{n_k})-1=n_k-1\ge k$. But this means
that $f$ is unbounded on the $C$-compact set $A$. So the sequence
$\{f_{n}\}_{n\in \mathbb N}$ cannot have a cluster point in
$C_{rc}(X)$ and consequently, $C_{rc}(X)$ fails to be a $q$-space.
Hence, X must be hemi-$C$-compact.

$(8)\Rightarrow(1)$. Here we need the well-known result which says
that if the topology of a locally convex Hausdorff space is
generated by a countable family of seminorms, then it is
metrizable. Now the locally convex topology on $C(X)$ generated by
the countable family of seminorms $\{p_{A_n} : n\in \mathbb N \}$
is metrizable and weaker than the $C$-compact-open topology.
However, since for each $C$-compact set $A$ in $X$, there exists
$A_{n}$ such that $A\subseteq A_{n}$, the locally convex topology
generated by the family of seminorms $\{p_{A} : A\in RC(X)\}$,
that is, the $C$-compact-open topology, is weaker than the
topology generated by the family of seminorms $\{p_{A_{n}} : n\in
\mathbb N \}$. Hence, $C_{rc}(X)$ is metrizable.

\end{proof}

\section{Separable and second countability }

\begin{theorem}\label{th15} For any space $X$ and $\lambda\in \{MC(X)$, $SC(X)$, $CC(X)$, $PS(X)$, $RC(X) \}$, the following are equivalent.

\begin{enumerate}

\item  $C_{p}(X)$ is separable.

\item  $C_{c}(X)$ is separable.

\item $X$ has a weaker separable metrizable topology.

\item  $C_{\lambda}(X)$ is separable.

\end{enumerate}

\end{theorem}

\begin{proof}
First, note by Corollary 4.2.2 in ~{\cite{mcnt}} that (1), (2),
and (3) are equivalent. Also, since $C_p(X)\le C_{\lambda}(X)$,
for $\lambda\in \{MC(X)$, $SC(X)$, $CC(X)$, $PS(X)$, $RC(X) \}$,
$(4)\Rightarrow(1)$.

$(3)\Rightarrow(4)$. If $X$ has a weaker separable metrizable
topology, then $X$ is submetrizable. By Theorem \ref{th10},
$C_{k}(X)=C_{c}(X)=C_{sc}(X)=C_{cc}(X)=C_{ps}(X)=C_{rc}(X)$.
 Since $(3)\Rightarrow(2)$, $C_{\lambda}(X)$ is separable for each $\lambda\in \{MC(X)$, $SC(X)$, $CC(X)$, $PS(X)$, $RC(X) \}$.

\end{proof}

\begin{corollary} If $X$ is pseudocompact and $\lambda\in \{MC(X)$, $K(X)$, $SC(X)$, $CC(X)$, $PS(X)$, $RC(X) \}$, then the following
statements are equivalent.

\begin{enumerate}

\item  $C_{\lambda}(X)$ is separable.

\item  $C_{\lambda}(X)$ has $ccc$.

\item  $X$ is metrizable.

\end{enumerate}

\end{corollary}

\begin{proof}
$(1)\Rightarrow(2)$. This is immediate.

$(2)\Rightarrow(3)$. By Corollary 4.8 in ~{\cite{os2}}, $X$ is
metrizable.

$(3)\Rightarrow(1)$. If $X$ is metrizable, then $X$, being
pseudocompact, is also compact. Hence $X$ is separable and
consequenly by Theorem \ref{th15}, $C_{\lambda}(X)$ is separable.

\end{proof}

Recall that a family of nonemty open sets in a space $X$ is called
a $\pi$-base for $X$ if every nonempty open set in $X$ contains a
member of this family.

The following  Theorems are analogues of Theorem 4.6 and Theorem
4.8  in ~{\cite{kuga1}}.

\begin{theorem} For a space $X$ and $\lambda\in \{MC(X)$, $K(X)$, $SC(X)$, $CC(X)$, $PS(X)$, $ RC(X) \}$, the following statements are
equivalent.

\begin{enumerate}

\item  $C_{\lambda}(X)$ contains a dense subspace which has a
countable $\pi$-base.

\item  $C_{\lambda}(X)$ has a countable $\pi$-base.

\item  $C_{\lambda}(X)$ is second countable.

\item  $X$ is hemicompact and $\aleph_0$-space.

\end{enumerate}

\end{theorem}

\begin{theorem} For a locally compact space $X$ and $\lambda\in \{MC(X)$, $K(X)$, $SC(X)$, $CC(X)$, $PS(X)$, $RC(X) \}$, the following
statements are equivalent.

\begin{enumerate}

\item  $C_{\lambda}(X)$ is second countable.

\item  $X$ is hemicompact and submetrizable.

\item  $X$ is Lindel\"{o}f and submetrizable.

\item  $X$ is the union of a countable family of compact
metrizable subsets of $X$.

\item  $X$ is second countable.

\end{enumerate}

\end{theorem}

\bibliographystyle{model1a-num-names}
\bibliography{<your-bib-database>}

\begin{thebibliography}{10}

\bibitem{ardug}
 R. Arens, J. Dugundji, \textit{Topologies for
 function spaces},  Pacific. J. Math.{\bf 1}, (1951), 5--31.

\bibitem{enge}
R. Engelking, \textit{General Topology}, PWN, Warsaw, (1977); Mir,
Moscow, (1986).



\bibitem{kund}
S. Kundu, A.B. Raha, \textit{The bounded-open topology and its
relatives}, Rend. Istit. Mat. Univ. Trieste {\bf 27} (1995),
61-77.

\bibitem{kuga}
S. Kundu, P. Garg, \textit{ The pseudocompact-open topology on
$C(X)$}, Topology Proceedings, VOL.~30, (2006), 279-299.

\bibitem{kuga1}
S. Kundu, P. Garg, \textit{ Countability properties of the
pseudocompact-open topology on $C(X)$: a comparative study}, Rend.
Istit. Mat. Univ. Trieste {\bf 39} (2007), 421-444.




\bibitem{mcnt}
 R.A. McCoy, I. Ntantu, \textit{Topological Properties of Spaces of Continuous
 Functions}, Lecture Notes in Math., 1315,
 Springer-Verlag, Berlin, (1988).

\bibitem{os2}
A.V. Osipov, \textit{ Topological-algebraic properties of function
spaces with set-open topologies}, Topology and its Applications,
{\bf 159}, issue 3, (2012), 800-805.


\bibitem{os1}
A.V. Osipov, \textit{ The Set-Open topology}, Top. Proc. {\bf 37}
(2011), 205-217.


\end{thebibliography}







\end{document}